\documentclass{article}
\usepackage{graphicx}
\usepackage{graphics, graphicx, latexsym, amssymb}
\newtheorem{propo}{Proposition}
\newtheorem{coro}{Corollary}

\begin{document}

\title{Polyominoes with nearly convex columns: \\
A semi-directed model}
\author{Svjetlan Fereti\'{c} \footnote{e-mail: svjetlan.feretic@gradri.hr} \\ 
Faculty of Civil Engineering, University of Rijeka, \\ 
Viktora Cara Emina 5, 51000 Rijeka, Croatia}
\maketitle

\begin{abstract}
Column-convex polyominoes are by now a well-explored model. So far, however, no attention has been given to polyominoes whose columns can have either one or two connected components. This little known kind of polyominoes seems not to be manageable as a whole. To obtain solvable models, one needs to introduce some restrictions. This paper is focused on polyominoes with hexagonal cells. The restrictions just mentioned are semi-directedness and an upper bound on the size of the gap within a column. The solvable models so obtained have rational area generating functions, as column-convex polyominoes do. However, the growth constants of the new models are $4.1149\ldots$ and more, whereas the growth constant of column-convex polyominoes is $3.8631\ldots \:$.
\end{abstract}

\vspace{8mm}
\noindent \textit{Keywords:} polyomino; hexagonal-celled; nearly convex column; semi-directed; area generating function

\vspace{3mm}
\noindent \textit{AMS Classification:} 05B50 (polyominoes); 05A15 (exact enumeration problems, generating functions)

\vspace{3mm}
\noindent \textit{Suggested running head:} Polyominoes with nearly convex columns

\newpage

\section{Introduction}

The enumeration of polyominoes is a topic of great interest to chemists, physicists and combinatorialists alike. In chemical terms, any polyomino (with hexagonal cells) is a possible benzenoid hydrocarbon. In physics, determining the number of $n$-celled polyominoes is related to the study of two-dimensional percolation phenomena. In combinatorics, polyominoes are of interest in their own right because several polyomino models have good-looking exact solutions.

Let $a_n$ be the number of $n$-celled polyominoes. An exact formula for $a_n$ seems unlikely to ever be found. However, there exist notable results on the quantity 
$c=\lim_{n\rightarrow\infty} \sqrt[n]{a_n}$. (This quantity is called the \textit{growth constant}.) In the case of polyominoes with hexagonal cells, V\"{o}ge and Guttmann \cite{Voege} gave a rigorous proof that the growth constant $c$ exists and satisfies the inequality $4.8049\leq c\leq 5.9047\: $. In the same paper, it is estimated that $c=5.1831478(17)$.  

In this paper, we are not going to improve the above-stated bounds on $c$. Instead, we are going to revisit \textit{column-convex polyominoes}. The first to study this now-familiar model was Temperley \cite{Temperley} in 1956. Whether with square cells or with hexagonal cells, column-convex polyominoes have a rational area generating function. When cells are hexagons, the growth constant of column-convex polyominoes is $3.8631\ldots$ (Klarner \cite{Klarner}). This growth constant remained a record (among polyomino models having \textit{reasonably simple} exact solutions) until 1982, when Dhar, Phani and Barma \cite{Dhar} invented \textit{directed site animals}. Directed site animals with step-set $A=\{(1,\, 0),\: (0,\, 1)\}$ can be viewed as a family of polyominoes with square cells, while directed site animals with step-set $B=\{(1,\, 0),\: (0,\, 1),\: (1,\, 1)\}$ can be viewed as a family of polyominoes with hexagonal cells. If the step-set is $B$, the growth constant of directed site animals is exactly $4$. Incidentally, with either of the step-sets $A$ and $B$, the area generating function of directed site animals is algebraic (which is rather surprising) and satisfies a simple quadratic equation. Later on, in 2002, Bousquet-M\'{e}lou and Rechnitzer \cite{Rechnitzer} introduced \textit{stacked directed animals} and \textit{multi-directed animals}. Those two models substantially generalize directed site animals. Whether the cells are squares or hexagons, the area generating function of stacked directed animals is degree-two algebraic, and the area generating function of multi-directed animals is not D-finite. When cells are hexagons, the growth constant of stacked directed animals is exactly $4.5$, and the growth constant of multi-directed animals is $4.5878\ldots \:$.

Although descended from directed animals, multi-directed animals are also a superset of column-convex polyominoes. (To be precise, this holds when the cells are hexagons. It is not quite clear whether multi-directed animals with square cells are a superset of column-convex polyominoes with square cells.) Besides multi-directed animals, there exist two other generalizations of column-convex polyominoes. Those two generalizations are called $m$\textit{-convex polygons} \cite{m-convex} and \textit{prudent polygons} \cite{prudent}. So far, however, $m$-convex polygons and prudent polygons have not been enumerated by area; they have only been enumerated by perimeter. 

The aim of this paper is to define a model which (a) generalizes column-convex polyominoes, (b) possesses a reasonably simple area generating function, and (c) possesses a high growth constant. In view of the facts stated above, we shall have to compete with just one pre-existing model, namely with multi-directed animals. 

In this paper, we actually introduce a sequence of generalizations of \textit{hexa\-gonal-celled} column-convex polyominoes. Namely, we define \textit{level} $m$ \textit{cheesy polyominoes} ($m=1,\: 2,\: 3,\ldots$). At every level, our new model has a rational area generating function. However, as level increases, those rational generating functions quickly gain in size. In computing generating functions, we go up to level $3$. The computations are done by using Bousquet-M\'{e}lou's \cite{Bousquet} and Svrtan's~\cite{Svrtan} ``turbo" version of the Temperley method \cite{Temperley}. The growth constants of level $1$, level $2$ and level $3$ cheesy polyominoes turn out to be $4.1149\ldots \:$, $4.2318\ldots \:$ and $4.2886\ldots$, respectively.

Thus, the growth constants of cheesy polyominoes are not as high as $4.5878\ldots$, the growth constant of multi-directed animals. However, we expect to obtain greater growth constants ($4.5$ or more) in our next two papers. Namely, we have in view two generalizations of level $m$ cheesy polyominoes. The names of those two generalizations are \textit{level} $m$ \textit{polyominoes with cheesy blocks} and \textit{level} $m$ \textit{column-subconvex polyominoes}. Incidentally, for every $m \in \mathbb{N}$, the area generating function of level $m$ polyominoes with cheesy blocks is rational, whereas the area generating function of level $m$ column-subconvex polyominoes is probably not D-finite. We are also planning to study similar generalizations of square-celled column-convex polyominoes. Let us mention, however, that our new hexagonal-celled models behave somewhat better than their square-celled counterparts. Solving a square-celled level $1$ model requires almost as much effort as solving a hexagonal-celled level $2$ model.

\section{Definitions and conventions}

There are three regular tilings of the Euclidean plane, namely the triangular tiling, the square tiling, and the hexagonal tiling. We adopt the convention that every square or hexagonal tile has two horizontal edges. In a regular tiling, a tile is often referred to as a \textit{cell}. A plane figure $P$ is a \textit{polyomino} if $P$ is a union of finitely many cells and the interior of $P$ is connected. See Figure 1. Observe that, if a union of \textit{hexagonal} cells is connected, then it possesses a connected interior as well.

\begin{figure}
\begin{center}
\includegraphics[width=55mm]{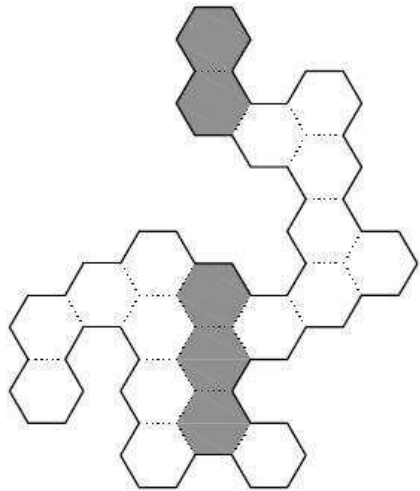}
\caption{A hexagonal-celled polyomino.}
\end{center}
\end{figure}

Let $P$ and $Q$ be two polyominoes. We consider $P$ and $Q$ to be equal if and only if there exists a translation $f$ such that $f(P)=Q$.

From now on, we concentrate on the hexagonal tiling. When we write ``a polyomino", we actually mean ``a hexagonal-celled polyomino".

Given a polyomino $P$, it is useful to partition the cells of $P$ according to their horizontal projection. Each block of that partition is a \textit{column} of $P$. Note that a column of a polyomino is not necessarily a connected set. An example of this is the highlighted column in Figure~1. On the other hand, it may happen that every column of a polyomino $P$ is a connected set. In this case, the polyomino $P$ is a \textit{column-convex polyomino}. See Figure 2.

\begin{figure}
\begin{center}
\includegraphics[width=47.5mm]{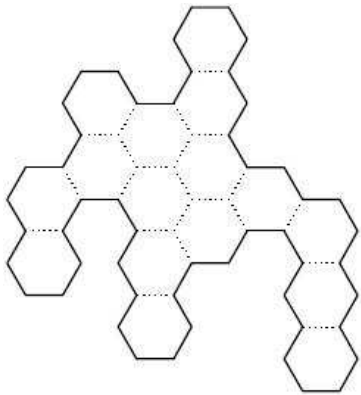}
\caption{A column-convex polyomino.}
\end{center}
\end{figure}

Let $a$ and $b$ be two adjacent columns of a polyomino $P$. Let the column $b$ be the right neighbour of the column $a$. It is certain that column $b$ has at least one edge in common with column $a$; otherwise $P$ could not be a connected set. However, there is no guarantee that every connected component of $b$ has at least one edge in common with $a$. For example, in Figure 1, the upper component of the highlighted column has no edge in common with the previous column.

Suppose that $P$ is such a polyomino that the first (\textit{i.e.}, leftmost) column of $P$ has no gap and that, in every pair of adjacent columns of $P$, every connected component of the right column has at least one edge in common with the left column. Then we say that $P$ is a \textit{rightward-semi-directed polyomino}.

A polyomino $P$ is a \textit{level} $m$ \textit{cheesy polyomino} if the following holds:

\begin{itemize}
\item $P$ is a rightward-semi-directed polyomino,
\item every column of $P$ has at most two connected components,
\item if a column of $P$ has two connected components, then the gap between the components consists of at most $m$ cells.
\end{itemize}

See Figure 3. 

\begin{figure}
\begin{center}
\includegraphics[width=55mm]{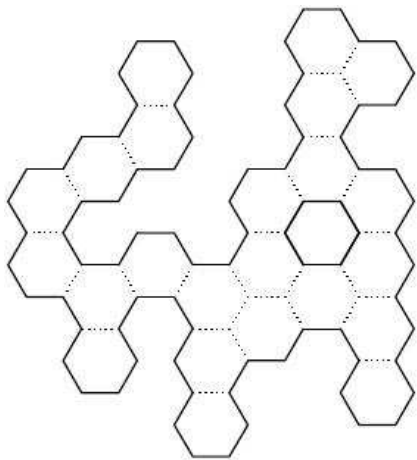}
\caption{A level one cheesy polyomino.}
\end{center}
\end{figure}

Level one cheesy polyominoes are a subset of level two cheesy polyominoes, level two cheesy polyominoes are a subset of level three cheesy polyominoes, and so on. As $m$ tends to infinity, the set of level $m$ cheesy polyominoes tends (in a certain sense) to the set of all polyominoes which are rightward-semi-directed and are made up of columns with at most two connected components.

The name ``cheesy polyominoes" is intended to suggest that these polyominoes can have internal holes. At this point, it might be objected that there exists another model, called \textit{directed animals}, in which holes occur more freely than in cheesy polyominoes. A column of a directed animal can have any number of holes, and the sizes of those holes are not subject to any limitations. However, the name of our model has 11 years of tradition\footnote{The contents of this paper exist since 1999. In that year, I defined cheesy polyominoes, explored them, and presented them at Mathematical Colloquium in Osijek \cite{Osijek_1999}, as well as at MATH/CHEM/COMP Course \& Conference in Dubrovnik \cite{Dubrovnik_1999}.}. So, our model will retain the name ``cheesy polyominoes" despite the fact that directed animals are arguably ``cheesier".

If a polyomino $P$ is made up of $n$ cells, we say that the \textit{area} of $P$ is $n$.

Let $a$ be a column of a polyomino $P$. By the \textit{height} of $a$ we mean the number of those cells which make up $a$ plus the number of those (zero or more) cells which make up the gaps of $a$. For example, in Figure 1, the highlighted column has height $7$, and the next column to the left has height $4$. 

Let $R$ be a set of polyominoes. By the \textit{area generating function} of $R$ we mean the formal sum

\begin{displaymath}
\sum_{P \in R} q^{area\ of\ P}.
\end{displaymath}

By the \textit{area and last column generating function} of $R$ we mean the formal sum

\begin{displaymath}
\sum_{P \in R} q^{area\ of\ P} \cdot t^{the\ height\ of\ the\ last\ column\ of\ P}.
\end{displaymath}

\section{Cheesy polyominoes vs. multi-directed animals}

This section has been removed because arXiv warned us that the submission exceeds size limits. The removed section contained definitions and some pictures of the directed classes mentioned in the introduction. From those definitions and pictures, it was clear that cheesy polyominoes are very different from the said directed models.

\section{Column-convex polyominoes}

\setcounter{figure}{8}

The area generating function for column-convex polyominoes is known since 1967 (Klarner \cite{Klarner}). However, rederiving that formula here will add to the completeness of this paper.

Let $A=A(q,t)$ be the area and last column generating function for column-convex polyominoes. Let $A_1=A(q,1)$ and $B_1=\frac{\partial A}{\partial t}(q,1)$. 

Let $S$ be the set of all column-convex polyominoes.

We write $S_{\alpha}$ for the set of column-convex polyominoes which have only one column. For $P \in S \setminus S_{\alpha}$, we define the \textit{pivot cell} of $P$ to be the lower right neighbour of the lowest cell of the second last column of $P$. See Figure~9. Observe that the pivot cell of a polyomino $P \in S \setminus S_{\alpha}$ is not necessarily contained in $P$. 

\begin{figure}
\begin{center}
\includegraphics[width=\textwidth]{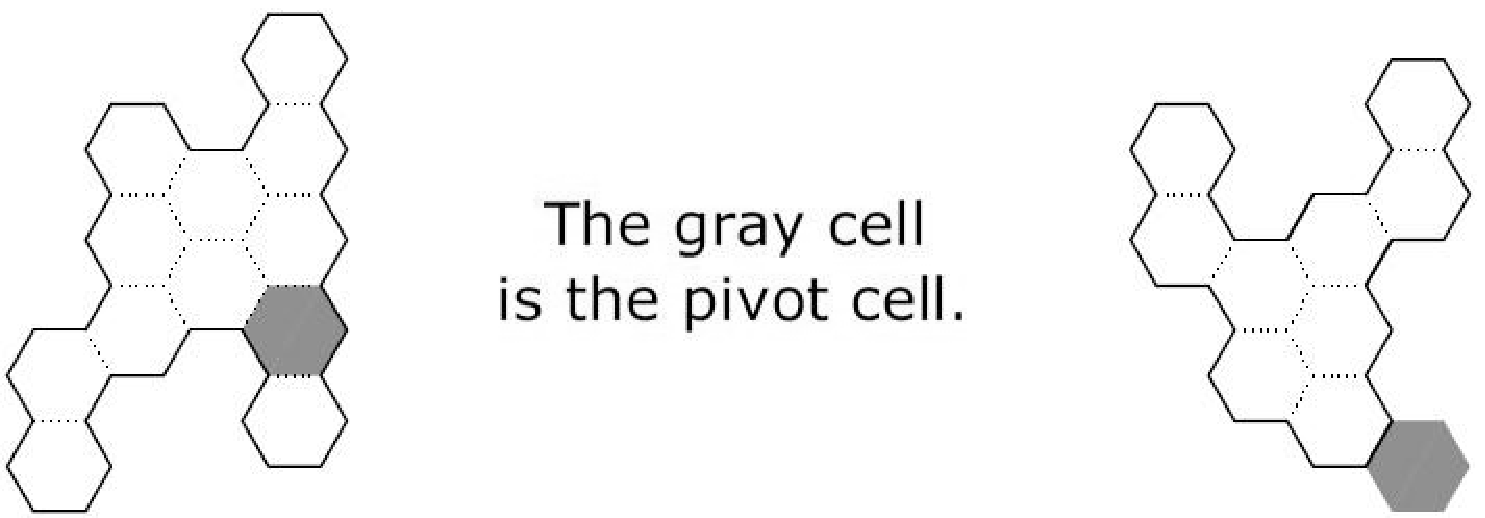}
\caption{The pivot cell.}
\end{center}
\end{figure}

Let

\begin{eqnarray*}
S_{\beta} & = & \{P \in S \setminus S_{\alpha}: P \mathrm{\ contains \ its \ pivot \ cell} \} \quad \mathrm{and} \\
S_{\gamma} & = & \{P \in S \setminus S_{\alpha}: P \mathrm{\ does \ not \ contain \ its \ pivot \ cell} \}.
\end{eqnarray*}

The sets $S_{\alpha}$, $S_{\beta}$ and $S_{\gamma}$ form a partition of $S$. We write $A_{\alpha}$, $A_{\beta}$ and $A_{\gamma}$ to denote the parts of the series $A$ that come from the sets $S_{\alpha}$, $S_{\beta}$ and $S_{\gamma}$, respectively.

It is obvious that

\begin{equation}
A_{\alpha}=qt+(qt)^2+(qt)^3+\ldots =\frac{qt}{1-qt}.
\end{equation}

If a polyomino $P$ lies in $S_{\beta}$, then the last column of $P$ is made up of the pivot cell, of $i \in \{0,\: 1,\: 2,\: 3,\ldots \: \}$ cells lying below the pivot cell, and of $j \in \{0,\: 1,\: 2,\: 3,\ldots \: \}$ cells lying above the pivot cell. Hence,

\begin{equation}
A_{\beta} = A_1 \cdot qt \cdot \left[ \sum_{i=0}^{\infty} (qt)^i \right] \cdot \left[ \sum_{j=0}^{\infty} (qt)^j \right] = \frac{qt}{(1-qt)^2} \cdot A_1.
\end{equation}

See Figure 10. 

\begin{figure}
\begin{center}
\includegraphics[width=40mm]{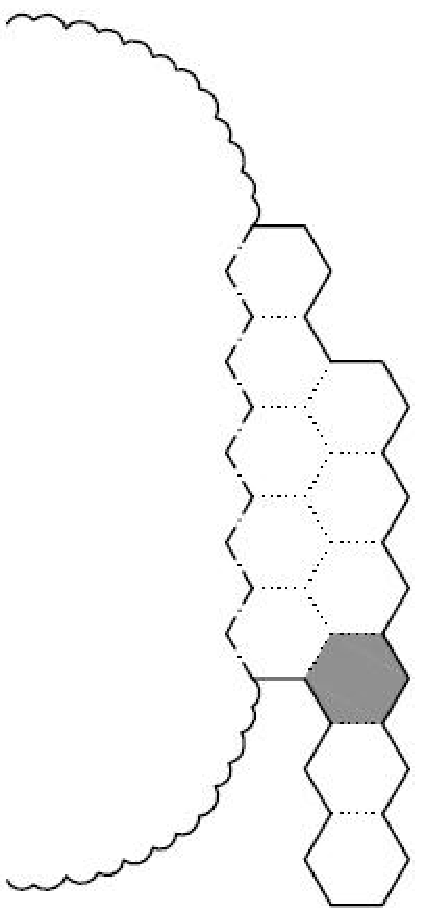}
\caption{The last two columns of an element of $S_{\beta}$.}
\end{center}
\end{figure}

For $n \in \{1,\: 2,\: 3,\ldots\}$, let $S^{(n)}$ denote the set of column-convex polyominoes whose last column consists of $n$ cells. Let $A_1^{(n)}$ be the part of $A_1$ that comes from the set $S^{(n)}$. Now, every element of $S_{\gamma}$ can be produced in three steps. Step one: We choose a number $n \in \{1,\: 2,\: 3,\ldots\}$ and a polyomino $P \in S^{(n)}$. Step two: In the last column of $P$, we choose a cell $c$. Step three: After choosing a number $i \in \{1,\: 2,\: 3,\ldots\}$, we place a new column of height $i$ so that the lowest cell of the new column is the upper right neighbour of the cell $c$. See Figure 11. Thus,

\begin{eqnarray}
& & A_{\gamma} = \sum_{n=1}^{\infty} A_1^{(n)} \cdot n \cdot \sum_{i=1}^{\infty} (qt)^i = \sum_{n=1}^{\infty} n \cdot A_1^{(n)} \cdot \frac{qt}{1-qt} \nonumber \\
& & = \frac{qt}{1-qt} \cdot \left[ A_1^{(1)}+2A_1^{(2)}+3A_1^{(3)}+\ldots \right] \nonumber \\ 
& & = \frac{qt}{1-qt} \cdot \left[ A_1^{(1)}+2A_1^{(2)}t+3A_1^{(3)}t^2+\ldots \right]_{with \ t=1} \nonumber \\
& & = \frac{qt}{1-qt} \cdot \left[ \frac{\partial A(q,t)}{\partial t} \right]_{with \ t=1} = \frac{qt}{1-qt} \cdot B_1.
\end{eqnarray}

\begin{figure}
\begin{center}
\includegraphics[width=66mm]{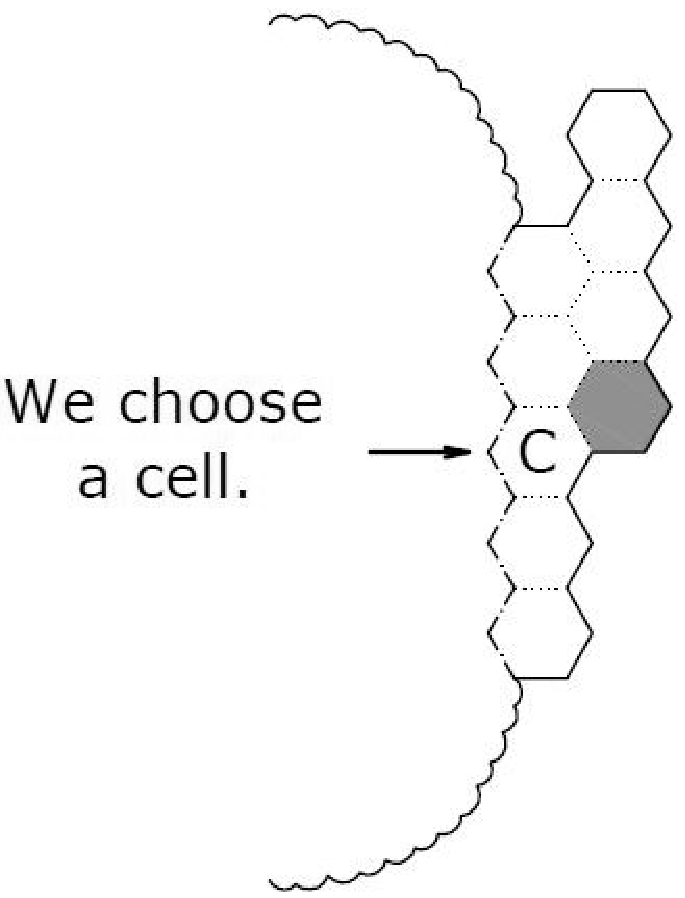}
\caption{The making of an element of $S_{\gamma}$.}
\end{center}
\end{figure}

Since $A=A_{\alpha}+A_{\beta}+A_{\gamma}$, Eqs. (1), (2) and (3) imply that

\begin{equation}
A=\frac{qt}{1-qt} + \frac{qt}{(1-qt)^2} \cdot A_1 + \frac{qt}{1-qt} \cdot B_1.
\end{equation}

Setting $t=1$, from Eq. (4) we obtain

\begin{displaymath}
A_1=\frac{q}{1-q} + \frac{q}{(1-q)^2} \cdot A_1 + \frac{q}{1-q} \cdot B_1.
\end{displaymath}

Differentiating Eq. (4) with respect to $t$ and then setting $t=1$, we obtain

\begin{displaymath}
B_1 = \frac{q}{1-q} + \frac{q^2}{(1-q)^2} + \left[ \frac{q}{(1-q)^2} + \frac{2q^2}{(1-q)^3} \right] \cdot A_1 + \left[ \frac{q}{1-q} + \frac{q^2}{(1-q)^2} \right] \cdot B_1.
\end{displaymath}

We now have a system of two linear equations in two unknowns, $A_1$ and $B_1$. By solving the system, we get the following proposition.

\begin{propo}
The area generating function for column-convex polyominoes is given by

\begin{displaymath}
A_1=\frac{q \cdot (1-q)^3}{1-6q+10q^2-7q^3+q^4} \ .
\end{displaymath}
\end{propo}

Being a quartic polynomial, the denominator of $A_1$ has four roots. Those roots are $r_1=0.2588\ldots$, $r_2=0.7066\ldots -0.4750\ldots \cdot i$, $r_3=0.7066\ldots +0.4750\ldots \cdot i$ and $r_4=5.3278\ldots \:$. The root with smallest absolute value is $r_1=0.2588\ldots$, and $\frac{1}{r_1}$ is equal to $3.8631\ldots \:$. Therefore, the coefficient of $q^n$ in $A_1$ (denoted $[q^n] A_1$) has the asymptotic behaviour $[q^n] A_1 \sim c \times 3.8631\ldots ^n$, where $c$ is a constant. To find the value of $c$, we decompose $A_1$ into partial fractions. The partial fraction involving $r_1$ turns out to be $-\frac{0.0487\ldots}{q-0.2588\ldots}$. In the Taylor series expansion of $-\frac{0.0487\ldots}{q-0.2588\ldots}$, the coefficient of $q^n$ is equal to $0.1884\ldots \times 3.8631\ldots ^n$. Thus, $c=0.1884\ldots \:$. We have got the following result.

\begin{coro} The number of $n$-celled column-convex polyominoes $[q^n] A_1$ has the asymptotic behaviour 
\begin{displaymath}
[q^n] A_1 \sim 0.188419\ldots \times 3.863130\ldots ^n.
\end{displaymath}
\end{coro}

Thus, the growth constant of column-convex polyominoes is $3.8631\ldots \:$.

\section{Level one cheesy polyominoes}

Counting level one cheesy polyominoes by area is nearly as easy as counting column-convex polyominoes by area.

Let $C=C(q,t)$ be the area and last column generating function for level one cheesy polyominoes. Let $C_1=C(q,1)$ and $D_1=\frac{\partial C}{\partial t}(q,1)$. 

Let $T$ be the set of all level one cheesy polyominoes. We write $T_{\alpha}$ for the set of level one cheesy polyominoes which have only one column. For $P \in T \setminus T_{\alpha}$, we define the \textit{pivot cell} of $P$ to be the lower right neighbour of the lowest cell of the second last column of $P$. As with column-convex polyominoes, the pivot cell of a polyomino $P \in T \setminus T_{\alpha}$ is not necessarily contained in $P$. Let

\begin{eqnarray*}
T_{\beta} & = & \{P \in T \setminus T_{\alpha}: \mathrm{the \ last \ column \ of \ } P \mathrm{\ has \ no \ hole,} \\ 
& & \mathrm{and \ the \ pivot \ cell \ of \ } P \mathrm{\ is \ contained \ in \ } P \}, \\
T_{\gamma} & = & \{P \in T \setminus T_{\alpha}: \mathrm{the \ last \ column \ of \ } P \mathrm{\ has \ no \ hole,} \\
& & \mathrm{and \ the \ pivot \ cell \ of \ } P \mathrm{\ is \ not \ contained \ in \ } P \}, \quad \mathrm{and} \\
T_{\delta} & = & \{P \in T \setminus T_{\alpha}: \mathrm{the \ last \ column \ of \ } P \mathrm{\ has \ a \ hole}\}. 
\end{eqnarray*}

The sets $T_{\alpha}$, $T_{\beta}$, $T_{\gamma}$ and $T_{\delta}$ form a partition of $T$. We write $C_{\alpha}$, $C_{\beta}$, $C_{\gamma}$ and $C_{\delta}$ to denote the parts of the series $C$ that come from the sets $T_{\alpha}$, $T_{\beta}$, $T_{\gamma}$ and $T_{\delta}$, respectively.

Similarly as with column-convex polyominoes, we have

\begin{equation}
C_{\alpha}=qt+(qt)^2+(qt)^3+\ldots =\frac{qt}{1-qt}.
\end{equation}

\noindent and

\begin{equation}
C_{\beta} = \frac{qt}{(1-qt)^2} \cdot C_1.
\end{equation}

Consider the following situation. A cheesy polyomino $P$ ends with a column $I$. We are creating a new column to the right of $I$, and the result should be an element of $T_{\gamma}$. Then, whether or not the column $I$ has a hole, we can put the lowest cell of the new column in exactly $m$ places, where $m$ is the height of $I$. See Figure 12. Hence

\begin{equation}
C_{\gamma}=\frac{qt}{1-qt} \cdot D_1.
\end{equation}

\begin{figure}
\begin{center}
\includegraphics[width=114mm]{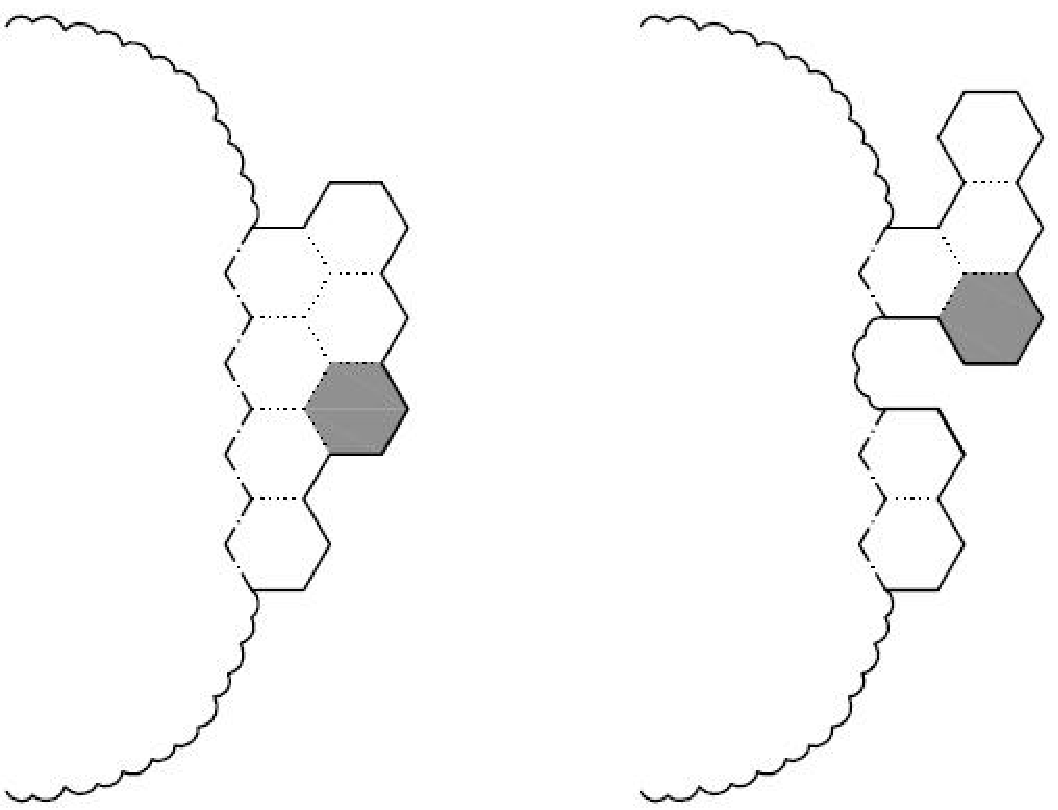}
\caption{The last two columns of two elements of $T_{\gamma}$.}
\end{center}
\end{figure}

Let us move on to another situation. A cheesy polyomino $P$ ends with a column $J$. We are creating a new column to the right of $J$, and the result should be an element of $T_{\delta}$. Then, whether or not the column $J$ has a hole, we can put the hole of the new column in exactly $n-1$ places, where $n$ is the height of $J$. See Figure 13. 

\begin{figure}
\begin{center}
\includegraphics[width=114mm]{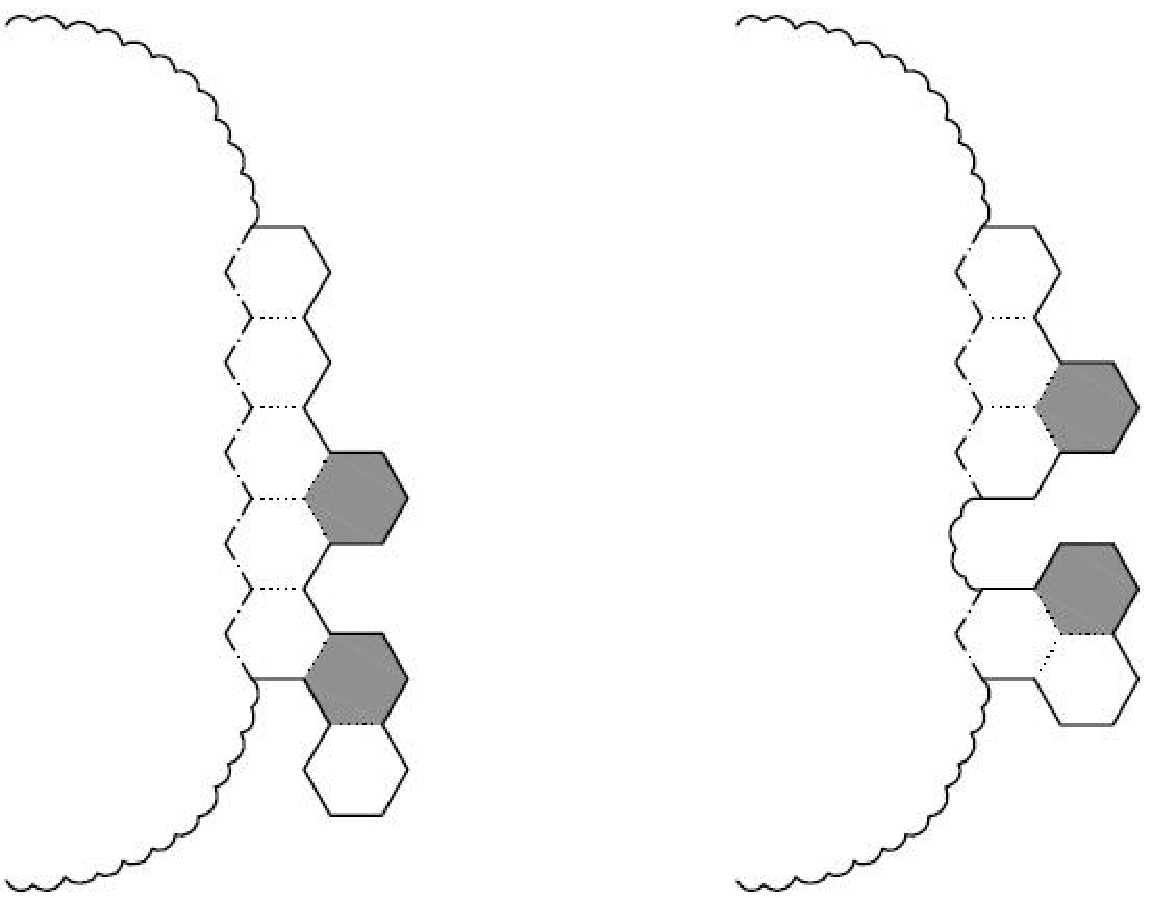}
\caption{The last two columns of two elements of $T_{\delta}$.}
\end{center}
\end{figure}

For $r \in \{1,\: 2,\: 3,\ldots\}$, let $T^{(r)}$ denote the set of level one cheesy polyominoes whose last column has height $r$. Let $C_1^{(r)}$ be the part of $C_1$ that comes from the set $T^{(r)}$. Every element of $T_{\delta}$ can be produced in four steps. Step one: We choose a number $r \in \{1,\: 2,\: 3,\ldots\}$ and a polyomino $P \in T^{(r)}$. Step two: In the last column of $P$, we choose two adjacent cells. (This can be done in $r-1$ ways.) Step three: We choose two numbers, $i$ and $j$, from the set $\{1,\: 2,\: 3,\ldots\}$. Step four: We create a two-component new column in which the upper component has $i$ cells, the lower component has $j$ cells, and the hole is the common right neighbour of the two cells chosen in Step two. Thus,

\begin{eqnarray}
& & C_{\delta} = \sum_{r=1}^{\infty} C_1^{(r)} \cdot (r-1) \cdot \left[ \sum_{i=1}^{\infty} (qt)^i \right] \cdot t \cdot \left[ \sum_{j=1}^{\infty} (qt)^j \right] \nonumber \\
& & = \sum_{r=1}^{\infty} (r-1) \cdot C_1^{(r)} \cdot \frac{q^2t^3}{(1-qt)^2} \nonumber \\
& & = \frac{q^2t^3}{(1-qt)^2} \cdot \left[ C_1^{(2)}+2C_1^{(3)}+3C_1^{(4)}+\ldots \right] \nonumber \\ 
& & = \frac{q^2t^3}{(1-qt)^2} \cdot \left\{ \left[C_1^{(1)}+2C_1^{(2)}t+3C_1^{(3)}t^2+\ldots \right]_{with \ t=1} \right. \nonumber \\
& & \left. \mbox{}-\left[C_1^{(1)}t+C_1^{(2)}t^2+C_1^{(3)}t^3+\ldots \right]_{with \ t=1} \right\} \nonumber \\
& & = \frac{q^2t^3}{(1-qt)^2} \cdot \left\{ \left[ \frac{\partial C(q,t)}{\partial t} \right]_{with \ t=1} - \left[ C(q,t)\right]_{with \ t=1} \right\} \nonumber \\
& & = \frac{q^2t^3}{(1-qt)^2} \cdot (D_1-C_1).
\end{eqnarray}

Since $C=C_{\alpha}+C_{\beta}+C_{\gamma}+C_{\delta}$, Eqs. (5), (6), (7) and (8) imply that

\begin{equation}
C=\frac{qt}{1-qt} + \frac{qt}{(1-qt)^2} \cdot C_1 + \frac{qt}{1-qt} \cdot D_1 + \frac{q^2 t^3}{(1-qt)^2} \cdot (D_1-C_1).
\end{equation}

Setting $t=1$, from Eq. (9) we obtain

\begin{displaymath}
C_1=\frac{q}{1-q} + \frac{q}{(1-q)^2} \cdot C_1 + \frac{q}{1-q} \cdot D_1 + \frac{q^2}{(1-q)^2} \cdot (D_1-C_1).
\end{displaymath}

Differentiating Eq. (9) with respect to $t$ and then setting $t=1$, we obtain

\begin{eqnarray*}
D_1 & = & \frac{q}{1-q} + \frac{q^2}{(1-q)^2} + \left[ \frac{q}{(1-q)^2} + \frac{2q^2}{(1-q)^3} \right] \cdot C_1 \\
& & \mbox{} + \left[ \frac{q}{1-q} + \frac{q^2}{(1-q)^2} \right] \cdot D_1 + \left[ \frac{3q^2}{(1-q)^2} + \frac{2q^3}{(1-q)^3} \right] \cdot (D_1-C_1).
\end{eqnarray*}

Thus, things are similar as with column-convex polyominoes. We have a system of two linear equations in two unknowns, $C_1$ and $D_1$. By solving the system, we get this proposition.

\begin{propo}
The area generating function for level one cheesy polyominoes is given by

\begin{displaymath}
C_1=\frac{q(1-3q+q^2)}{1-6q+8q^2-q^3}.
\end{displaymath}
\end{propo}

The roots of the denominator of $A_1$ are $r_1=0.2430\ldots$, $r_2=0.5727\ldots$ and $r_3=7.1842\ldots \:$. (All the three roots are real.) The root with smallest absolute value is $r_1=0.2430\ldots$, and $\frac{1}{r_1}$ is equal to $4.1149\ldots \:$. By decomposing $C_1$ into partial fractions and expanding the partial fractions into Taylor series, we find out the following fact.

\begin{coro}
The number of $n$-celled level one cheesy polyominoes $[q^n] C_1$ has the asymptotic behaviour

\begin{displaymath}
[q^n] C_1 \sim 0.144176\ldots \times 4.114907\ldots ^n.
\end{displaymath}
\end{coro}

Thus, the growth constant of level one cheesy polyominoes is $4.1149\ldots \:$. We observe a considerable increase with respect to $3.8631\ldots$, the growth constant of column-convex polyominoes.

\section{Level two cheesy polyominoes}

In this enumeration, if the last column of a polyomino has two connected components, we sometimes need to record not only the overall height of the last column, but also the height of the last column's upper component and the height of the last column's lower component. Hence, in addition to the ``old" variables $q$ and $t$, we introduce two new variables, $u$ and $v$. As before, the exponent of $q$ is the area and the exponent of $t$ is the overall height of the last column\footnote{Recall what do we mean by the height of a column: in Figure 1, the highlighted column has height $7$, and the next column to the left has height $4$.}. The exponent of $u$ is the height of the upper component of the last column, and the exponent of $v$ is the height of the lower component of the last column. 

The two main generating functions in this enumeration are $E=E(q,t)$ and $G=G(q,t,u,v)$. Those generating functions are used for the following purposes: 

\begin{itemize}
\item $E$ is a generating function for level two cheesy polyominoes whose last column either has no hole or has a one-celled hole,
\item $G$ is a generating function for level two cheesy polyominoes whose last column has a two-celled hole.
\end{itemize}

Let $E_1=E(q,1)$, $F_0=\frac{\partial E}{\partial t}(q,0)$, $F_1=\frac{\partial E}{\partial t}(q,1)$, $G_1=G(q,1,1,1)$, $H_1=\frac{\partial G}{\partial t}(q,1,1,1)$,
$I_0=\frac{\partial G}{\partial u}(q,1,0,1)$, and $J_0=\frac{\partial G}{\partial v}(q,1,1,0)$.

Let $U$ be the set of those level two cheesy polyominoes whose last column either has no hole or has a one-celled hole. Let $V$ be the set of those level two cheesy polyominoes whose last column has a two-celled hole. For $P \in U \cup V$, we define the \textit{body} of $P$ to be all of $P$, except the rightmost column of $P$. 

We are now going to partition the set $U$ into seven subsets and the set $V$ into two subsets. Let

\begin{eqnarray*}
U_{\alpha} & = & \{P \in U: P \mathrm{\ has\ only\ one\ column} \}, \\
U_{\beta} & = & \{P \in U \setminus U_{\alpha}: \mathrm{the\ body\ of\ } P \mathrm{\ lies\ in\ } U \mathrm{,\ the\ last\ column\ of\ } P \mathrm{\ has\ no\ hole,} \\ & & \mathrm{and\ the\ pivot\ cell\ of\ } P \mathrm{\ is\ contained\ in\ } P \}, \\
U_{\gamma} & = & \{P \in U \setminus U_{\alpha}: \mathrm{the\ body\ of\ } P \mathrm{\ lies\ in\ } U \mathrm{,\ the\ last\ column\ of\ } P \mathrm{\ has\ no\ hole,} \\ 
& & \mathrm{and\ the\ pivot\ cell\ of\ } P \mathrm{\ is\ not\ contained\ in\ } P \}, \\
U_{\delta} & = & \{P \in U \setminus U_{\alpha}: \mathrm{the\ body\ of\ } P \mathrm{\ lies\ in\ } U \mathrm{,\ and\ the\ last\ column\ of\ } P 
\mathrm{\ has\ a\ hole} \}, \\
U_{\epsilon} & = & \{P \in U \setminus U_{\alpha}: \mathrm{the\ body\ of\ } P \mathrm{\ lies\ in\ } V \mathrm{,\ the\ last\ column\ of\ } P \mathrm{\ has\ no\ hole,} \\ 
& & \mathrm{and\ the\ pivot\ cell\ of\ } P \mathrm{\ is\ contained\ in\ } P \}, \\
U_{\zeta} & = & \{P \in U \setminus U_{\alpha}: \mathrm{the\ body\ of\ } P \mathrm{\ lies\ in\ } V \mathrm{,\ the\ last\ column\ of\ } P \mathrm{\ has\ no\ hole,} \\ & & \mathrm{and\ the\ pivot\ cell\ of\ } P \mathrm{\ is\ not\ contained\ in\ } P \}, \quad \mathrm{and} \\
U_{\eta} & = & \{P \in U \setminus U_{\alpha}: \mathrm{the\ body\ of\ } P \mathrm{\ lies\ in\ } V \mathrm{,\ and\ the\ last\ column\ of\ } P 
\mathrm{\ has\ a\ hole} \}.
\end{eqnarray*}

Let

\begin{eqnarray*}
V_{\alpha} & = & \{P \in V: \mathrm{the \ body \ of \ } P \mathrm{\ lies \ in \ } U \}, \quad \mathrm{and} \\
V_{\beta} & = & \{P \in V: \mathrm{the \ body \ of \ } P \mathrm{\ lies \ in \ } V \}.
\end{eqnarray*}

It is clear that the sets $U_{\alpha},\ U_{\beta},\ldots,\ U_{\eta}$ form a partition of $U$, and that the sets $V_{\alpha}$ and $V_{\beta}$ form a partition of $V$. We shall write $E_{\alpha},\ E_{\beta},\ldots,\ E_{\eta}$ for the parts of the series $E$ that come from the sets $U_{\alpha},\ U_{\beta},\ldots,\ U_{\eta}$, respectively. Also, we shall write $G_{\alpha}$ and $G_{\beta}$ for the parts of the series $G$ that come from the sets $V_{\alpha}$ and $V_{\beta}$, respectively.

As in Section 5, we have

\begin{eqnarray}
E_{\alpha} & = & \frac{qt}{1-qt} \ , \\
E_{\beta} & = & \frac{qt}{(1-qt)^2} \cdot E_1, \\
E_{\gamma} & = & \frac{qt}{1-qt} \cdot F_1, \\
E_{\delta} & = & \frac{q^2t^3}{(1-qt)^2} \cdot (F_1-E_1), \\
E_{\epsilon} & = & \frac{qt}{(1-qt)^2} \cdot G_1.
\end{eqnarray}

The functional equation for $E_{\zeta}$ is more interesting. Let $P$ be an element of $V$, let $I$ be the last column of $P$, and let $m$ be the height of $I$. Suppose that we are creating a new column to the right of $I$, and that the result should be an element of $U_{\zeta}$. To the right of $I$ there is one cell (say $c$) which shares an edge with each of the two cells forming the hole of $I$. If we choose $c$ as the lowest cell of the new column, then the new column will have to be at least two cells high. Otherwise $P \cup \mathrm{(the\ new\ column)}$ will not be a polyomino. See Figure 14. In addition to $c$, there are $m-1$ other choices for the lowest cell of the new column. For each of these $m-1$ choices, $P \cup \mathrm{(the\ new\ column)}$ is a polyomino regardless of how many cells the new column has. Altogether, we have 

\begin{equation}
E_{\zeta} = \frac{q^2t^2}{1-qt} \cdot G_1 + \frac{qt}{1-qt} \cdot (H_1-G_1).
\end{equation}

\begin{figure}
\begin{center}
\includegraphics[width=114mm]{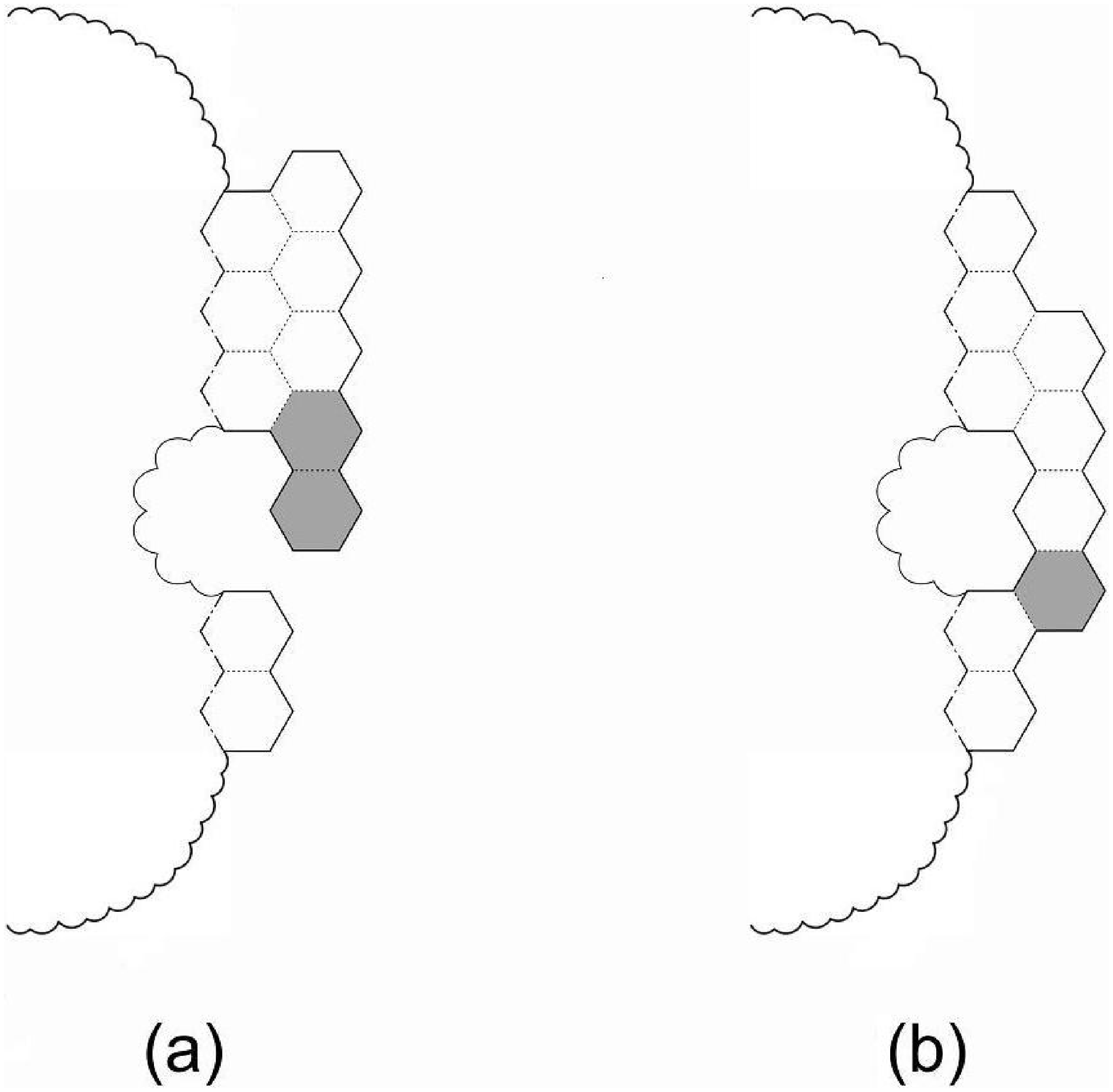}
\caption{(a) For this choice of the lowest cell, the new column must be at least two cells high. (b) For this choice of the lowest cell, the new column does not have to be at least two cells high.}
\end{center}
\end{figure}

We proceed to $E_{\eta}$. Let $P$ be an element of $V$, let $J$ be the last column of $P$, and let $n$ be the height of $J$. Suppose that we are creating a new column to the right of $J$, and that the result should be an element of $U_{\eta}$. One of the choices for the hole of the new column is the upper right neighbour of the top cell of the lower component of $J$. For this choice, the upper component of the new column has to have at least two cells. Otherwise $P \cup \mathrm{(the\ new\ column)}$ is not a polyomino. See Figure 15. We can also choose the hole of the new column as the lower right neighbour of the bottom cell of the upper component of $J$. Then, in order for $P \cup \mathrm{(the\ new\ column)}$ to be a polyomino, the lower component of the new column has to have at least two cells. In addition to the two ways just considered, there exist $n-3$ other ways to choose the hole of the new column. For each of those $n-3$ ways, $P \cup \mathrm{(the\ new\ column)}$ is a polyomino regardless of the sizes of the upper and lower components of the new column. Thus, we have 

\begin{equation}
E_{\eta} = \frac{2q^3t^4}{(1-qt)^2} \cdot G_1 + \frac{q^2t^3}{(1-qt)^2} \cdot (H_1-3G_1).
\end{equation}

\begin{figure}
\begin{center}
\includegraphics[width=114mm]{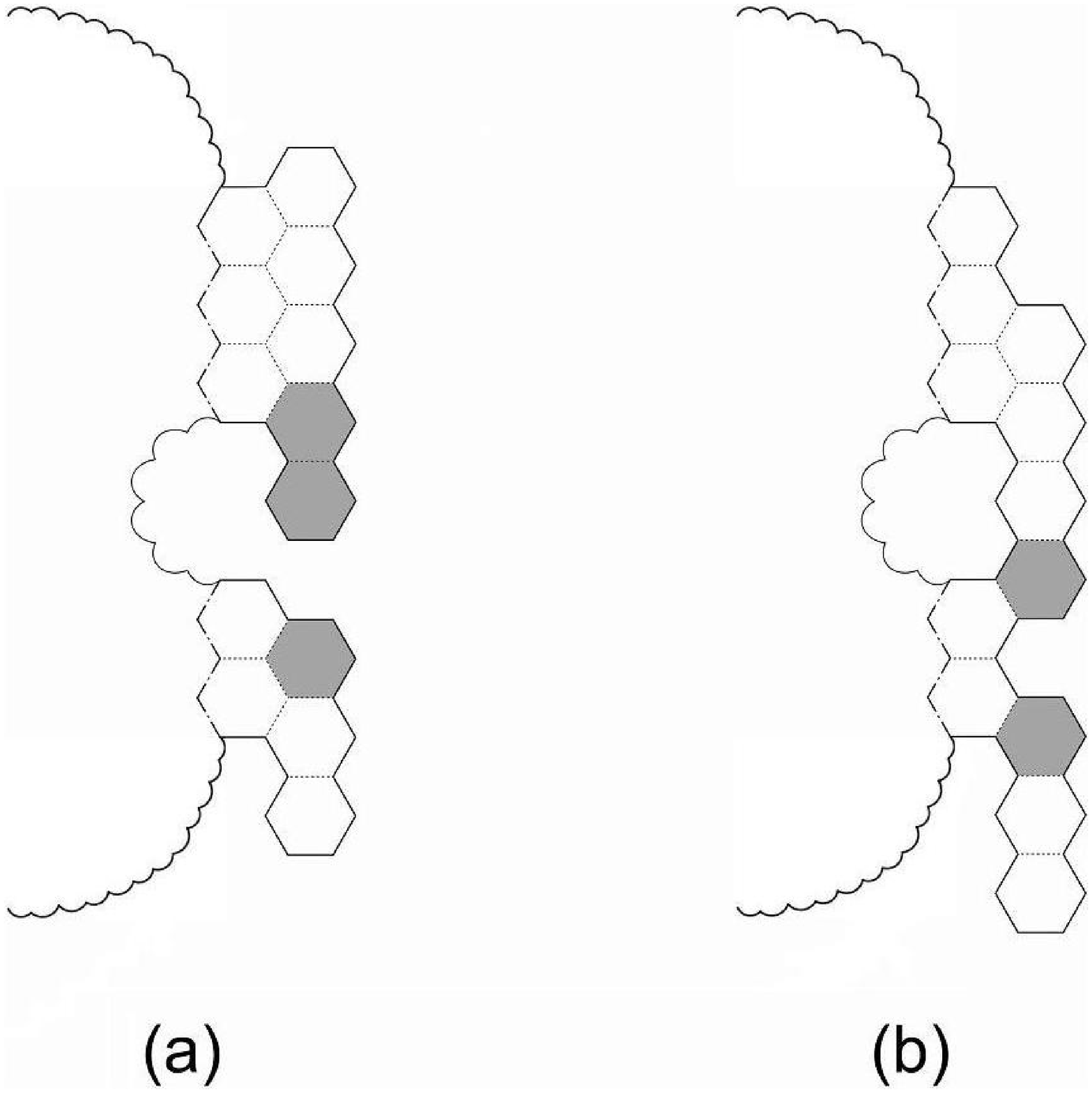}
\caption{(a) If the one-celled hole is in this position, the upper component of the new column must have at least two cells. (b) If the one-celled hole is in this position, the upper and lower components of the new column can be of any sizes.}
\end{center}
\end{figure}

Of course, the series $E_{\alpha},\ E_{\beta},\ldots,\ E_{\eta}$ sum up to $E$. Therefore, the summation of Eqs. (10)--(16) gives

\begin{eqnarray}
E & = & \frac{qt}{1-qt} + \frac{qt}{(1-qt)^2} \cdot E_1 + \frac{qt}{1-qt} \cdot F_1 + \frac{q^2t^3}{(1-qt)^2} \cdot (F_1-E_1) \nonumber \\
& & \mbox{} + \frac{qt}{(1-qt)^2} \cdot G_1 + \frac{q^2t^2}{1-qt} \cdot G_1 + \frac{qt}{1-qt} \cdot (H_1-G_1) \nonumber \\
& & \mbox{} + \frac{2q^3t^4}{(1-qt)^2} \cdot G_1 + \frac{q^2t^3}{(1-qt)^2} \cdot (H_1-3G_1).
\end{eqnarray}

We also need to establish a functional equation for $G$. Let $P$ be an element of $U$ and let $c$ be a column with a two-celled hole. Suppose that we want to glue $c$ to $P$ so that $P \cup c$ lies in $V_{\alpha}$, and so that $P$ and $c$ are the body and the last column of $P \cup c$, respectively. In how many ways $P$ and $c$ can be glued together? In principle, the number of ways is (the height of the last column of $P$) minus two. See Figure 16. However, if $P$ ends with a one-celled column, then we can glue $c$ to $P$ in zero ways, and not in minus one ways. Thus, we have

\begin{equation}
G_{\alpha}=\frac{q^2t^4uv}{(1-qtu)(1-qtv)} \cdot (F_1-2E_1+F_0).
\end{equation}

\begin{figure}
\begin{center}
\includegraphics[width=114mm]{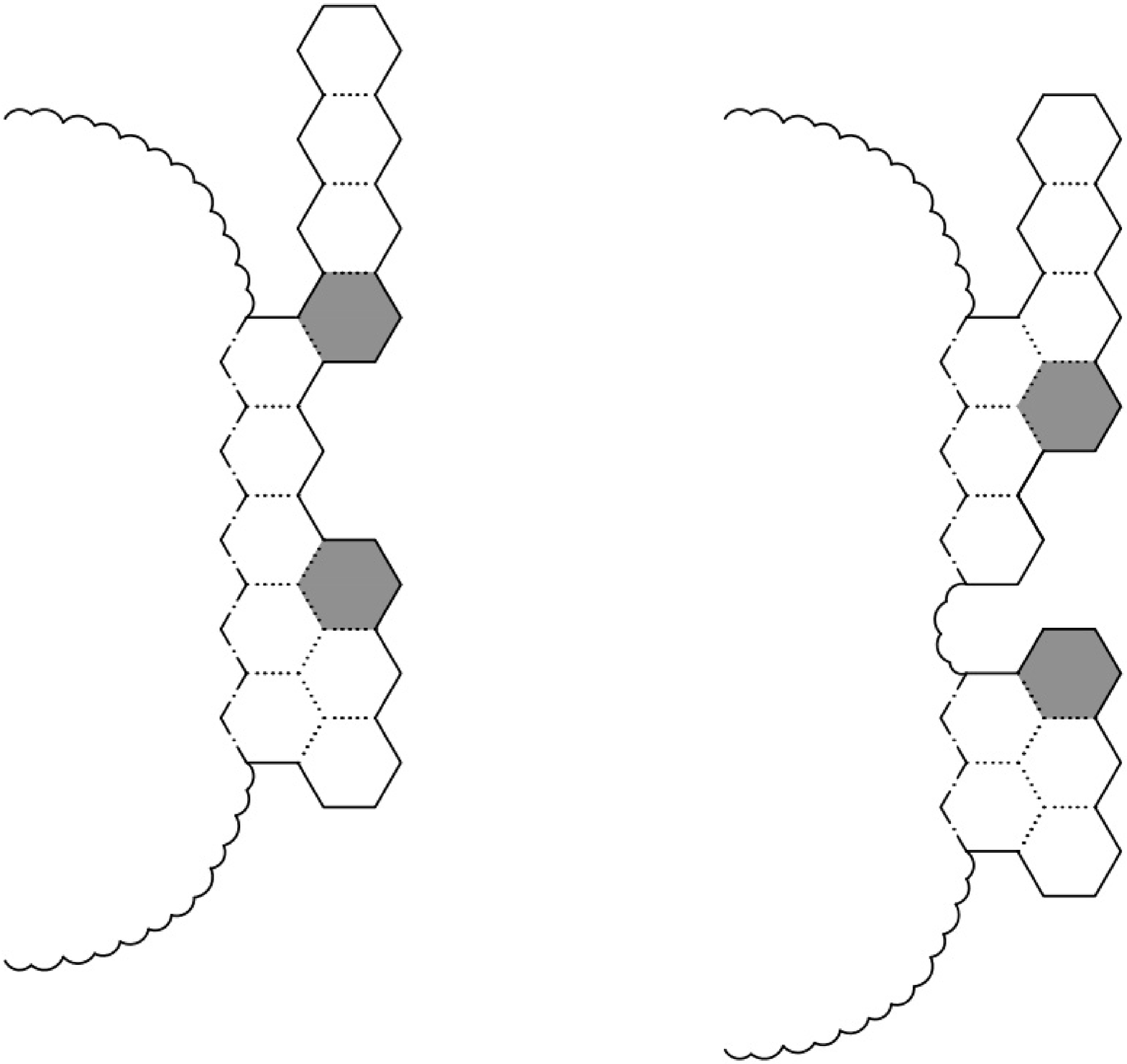}  
\caption{The last two columns of two elements of $V_{\alpha}$.}
\end{center}
\end{figure}

In the case of $G_{\beta}$, it is convenient to use overcounting. That is, we are going to ``mistakenly" count too much, and then subtract the parts which do not belong. Let $P$ be an element of $V$ and let $c$ be a column with a two-celled hole. Suppose that we want to glue $c$ to $P$ so that $P \cup c$ lies in $V_{\beta}$, and so that $P$ and $c$ are the body and the last column of $P \cup c$, respectively. In how many ways $P$ and $c$ can be glued together? First, there are (the height of the last column of $P$) minus two ways to satisfy these two necessary conditions:

\begin{itemize}
\item the bottom cell of the upper component of $c$ is either identical with or lies lower than the upper right neighbour of the top cell of the last column of $P$, and
\item the top cell of the lower component of $c$ is either identical with or lies higher than the lower right neighbour of the bottom cell of the last column of $P$. 
\end{itemize}

See Figure 17. 

\begin{figure}
\begin{center}
\includegraphics[width=118mm]{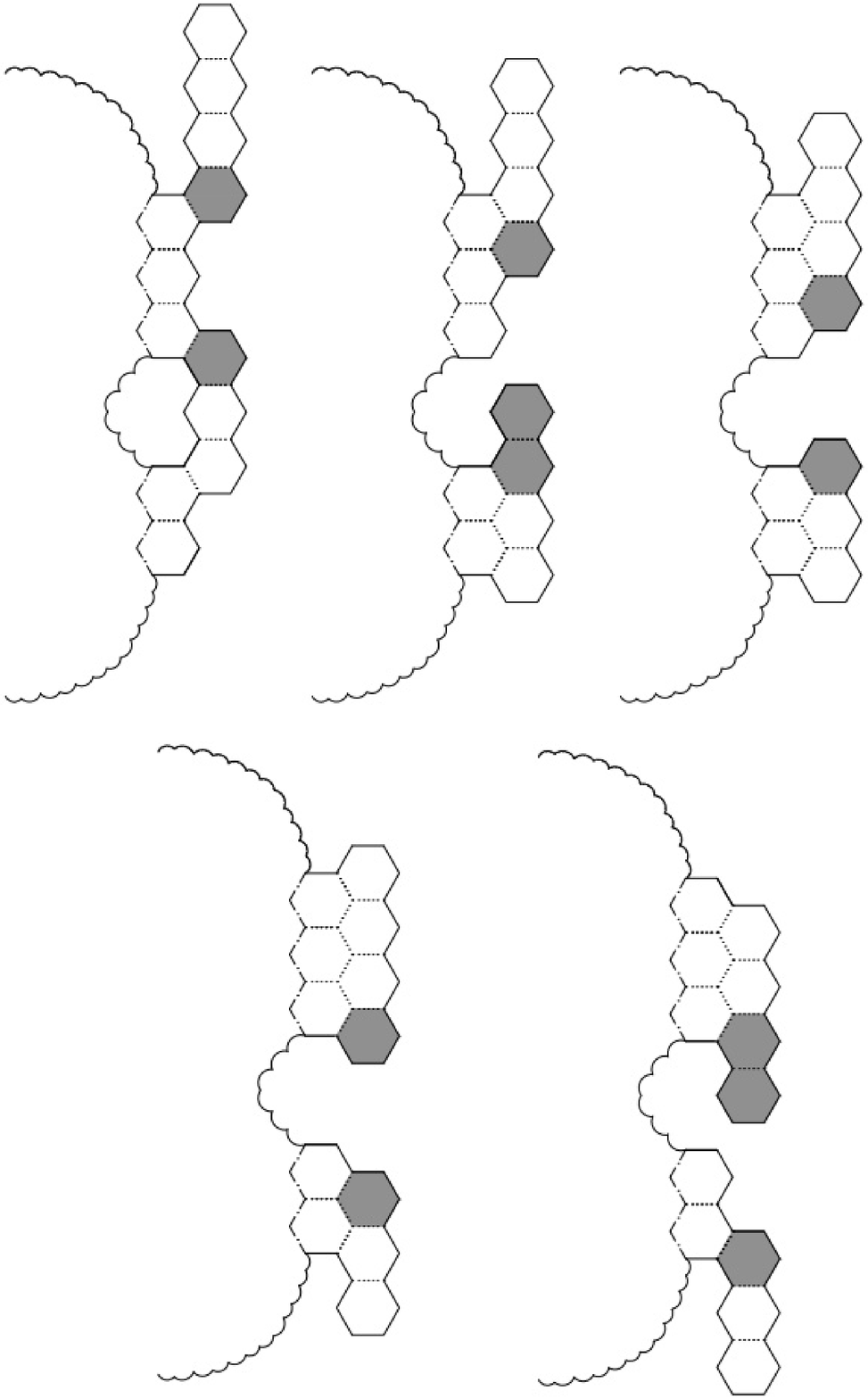}  
\caption{The last two columns of five elements of $V_{\beta}$.}
\end{center}
\end{figure}

So, if there were no special cases, then $G_{\beta}$ would be equal to

\begin{equation}
\frac{q^2t^4uv}{(1-qtu)(1-qtv)} \cdot (H_1-2G_1).
\end{equation}

However, special cases do exist. There are two of them:

\begin{enumerate}
\item The upper component of the last column of $P \in V$ has at least two cells and the lower component of the two-component column $c$ has just one cell.
\item The lower component of the last column of $P \in V$ has at least two cells and the upper component of the two-component column $c$ has just one cell.
\end{enumerate}

In case 1, it is (so to speak) dangerous to glue $c$ to $P$ in such a way that the one-celled lower component of $c$ becomes a common neighbour of the two cells which form the hole of the last column of $P$. This dangerous operation produces an object which is not a polyomino and hence does not lie in $V_{\beta}$.

In case 2, it is dangerous to glue $c$ to $P$ in such a way that the one-celled upper component of $c$ becomes a common neighbour of the two cells which form the hole of the last column of $P$. Again, the dangerous operation produces an object which is not a polyomino and hence does not lie in $V_{\beta}$.

Now, Eq. (19) is actually a generating function for the union of $V_{\beta}$ with the set of objects produced by the two dangerous operations. The generating function for the objects produced by the first dangerous operation is

\begin{equation}
\frac{q^2t^4uv}{1-qtu} \cdot (G_1-I_0).
\end{equation}

The generating function for the objects produced by the second dangerous operation is

\begin{equation}
\frac{q^2t^4uv}{1-qtv} \cdot (G_1-J_0).
\end{equation}

Subtracting Eqs. (20) and (21) from Eq. (19), we obtain

\begin{equation}
G_{\beta}=\frac{q^2t^4uv}{(1-qtu)(1-qtv)} \cdot (H_1-2G_1) - \frac{q^2t^4uv}{1-qtu} \cdot (G_1-I_0) - \frac{q^2t^4uv}{1-qtv} \cdot (G_1-J_0).
\end{equation}

Since $G=G_{\alpha}+G_{\beta}$, Eqs. (18) and (22) imply that

\begin{eqnarray}
G & = & \frac{q^2t^4uv}{(1-qtu)(1-qtv)} \cdot (F_1-2E_1+F_0+H_1-2G_1) \nonumber \\
& & \mbox{}-\frac{q^2t^4uv}{1-qtu} \cdot (G_1-I_0) - \frac{q^2t^4uv}{1-qtv} \cdot (G_1-J_0).
\end{eqnarray}

Using the computer algebra system \textit{Maple}, from Eqs. (17) and (23) we obtained a system of seven linear equations in seven unknowns: $E_1$, $F_0$, $F_1$, $G_1$, $H_1$, $I_0$ and $J_0$. Rather than write down these seven equations (some of which are a bit cumbersome), here below we give a list of recipes. Recipe no. $k$ ($k=1,\ 2,\ldots,\ 7$) tells how to obtain the $k$th equation of the linear system.

\begin{enumerate}
\item In Eq. (17), set $t=1$.
\item Differentiate Eq. (17) with respect to $t$ and then set $t=0$.
\item Differentiate Eq. (17) with respect to $t$ and then set $t=1$.
\item In Eq. (23), set $t=u=v=1$.
\item Differentiate Eq. (23) with respect to $t$ and then set $t=u=v=1$.
\item Differentiate Eq. (23) with respect to $u$. Then set $t=v=1$ and $u=0$.
\item Differentiate Eq. (23) with respect to $v$. Then set $t=u=1$ and $v=0$.
\end{enumerate}

The computer algebra quickly solved the linear system and then summed the generating functions $E_1$ and $G_1$. The result can be seen in the following proposition.

\begin{propo}
The area generating function for level two cheesy polyominoes is given by

\begin{displaymath}
K=\frac{q\cdot(1-6q+11q^{2}-6q^{3}-q^{4}-3q^{6}+5q^{7}+4q^{8}-3q^{9}-3q^{10})}{1-9q+27q^{2}-31q^{3}+8q^{4}+4q^{5}-2q^{6}+16q^{7}-5q^{8}-16q^{9}-2q^{10}+5q^{11}} \: .
\end{displaymath}
\end{propo}

\begin{coro}
The number of $n$-celled level two cheesy polyominoes has the asymptotic behaviour

\begin{displaymath}
[q^n] K \sim 0.121042\ldots \times 4.231836\ldots ^n.
\end{displaymath}
\end{coro}

Thus, the growth constant of level two cheesy polyominoes is $4.2318\ldots \:$.

\section{Level three cheesy polyominoes}

The enumeration of level three cheesy polyominoes is similar to the enumeration of level two cheesy polyominoes. However, there are still more cases than before. In the enumeration of level two cheesy polyominoes, the total number of cases was $9$ (because we partitioned the set $U$ into $7$ subsets and the set $V$ into $2$ subsets). At level three, the total number of cases is $16$. We deem it reasonable to skip those $16$ cases and only state the final result.

\begin{propo}
The area generating function for level three cheesy polyominoes is given by

\begin{displaymath}
L=\frac{M}{N} \: ,
\end{displaymath}

\noindent where

\begin{eqnarray*}
M & = & q \cdot (1-11q+49q^{2}-114q^{3}+146q^{4}-94q^{5}+5q^{6}+71q^{7}-143q^{8} \\
& & \mbox{}+176q^{9}-154q^{10}+100q^{11}+24q^{12}-121q^{13}+90q^{14}-61q^{15}+19q^{16} \\
& & \mbox{}+58q^{17}-32q^{18}-31q^{19}+37q^{20}+14q^{21}-43q^{22}-4q^{23}+21q^{24} \\
& & \mbox{}-q^{25}-5q^{26}) \\
& & \\
and & & \\
& & \\
N & = & 1-14q+80q^{2}-243q^{3}+423q^{4}-413q^{5}+174q^{6}+106q^{7}-350q^{8} \\
& & \mbox{}+533q^{9}-546q^{10}+427q^{11}-148q^{12}-261q^{13}+383q^{14}-253q^{15} \\
& & \mbox{}+158q^{16}+57q^{17}-181q^{18}+10q^{19}+115q^{20}-49q^{21}-96q^{22} \\
& & \mbox{}+93q^{23}+49q^{24}-54q^{25}-12q^{26}+12q^{27}+q^{28}.
\end{eqnarray*}
\end{propo}

\begin{coro}
The number of $n$-celled level three cheesy polyominoes has the asymptotic behaviour

\begin{displaymath}
[q^n] L \sim 0.108269\ldots \times 4.288630\ldots^n.
\end{displaymath}
\end{coro}

\section{Taylor expansions and the limit value of the growth constants}

To see how many polyominoes of a given type have $1,\: 2,\: 3,\ldots$ cells, we expanded the area generating functions into Taylor series. The results are shown in Table~1.

\begin{table}
\begin{center}
\begin{tabular}{|r||r|r|r|r|r|}\hline
& \multicolumn{1}{c|}{Column-} & \multicolumn{1}{c|}{Level 1} & \multicolumn{1}{c|}{Level 2} & \multicolumn{1}{c|}{Level 3} & \\
& \multicolumn{1}{c|}{convex} & \multicolumn{1}{c|}{cheesy} & \multicolumn{1}{c|}{cheesy} & \multicolumn{1}{c|}{cheesy} & \multicolumn{1}{c|}{All} \\
& \multicolumn{1}{c|}{polyo-} & \multicolumn{1}{c|}{polyo-} & \multicolumn{1}{c|}{polyo-} & \multicolumn{1}{c|}{polyo-} & \multicolumn{1}{c|}{polyo-} \\
\multicolumn{1}{|c||}{Area} & \multicolumn{1}{c|}{minoes} & \multicolumn{1}{c|}{minoes} & \multicolumn{1}{c|}{minoes} & \multicolumn{1}{c|}{minoes} & \multicolumn{1}{c|}{minoes} \\ \hline \hline
1 & 1 & 1 & 1 & 1 & 1 \\ \hline
2 & 3 & 3 & 3 & 3 & 3 \\ \hline
3 & 11 & 11 & 11 & 11 & 11 \\ \hline
4 & 42 & 43 & 43 & 43 & 44 \\ \hline
5 & 162 & 173 & 174 & 174 & 186 \\ \hline
6 & 626 & 705 & 718 & 719 & 814 \\ \hline
7 & 2419 & 2889 & 2996 & 3012 & 3652 \\ \hline
8 & 9346 & 11867 & 12579 & 12727 & 16689 \\ \hline
9 & 36106 & 48795 & 52996 & 54067 & 77359 \\ \hline
10 & 139483 & 200723 & 223705 & 230464 & 362671 \\ \hline
11 & 538841 & 825845 & 945324 & 984477 & 1716033 \\ \hline
12 & 2081612 & 3398081 & 3997185 & 4211222 & 8182213 \\ \hline
\end{tabular}
\caption{Here is how many polyominoes of a given type have $1,\: 2,\ldots,\: 12$ cells.}
\end{center}
\end{table}

The row ``$\mathrm{area}=4$" of Table 1 reads $42,\ 43,\ 43,\ 43,\ 44$. Indeed, Figure 18 shows the only two four-celled polyominoes which are not column-convex polyominoes. The polyomino on the left is a level $m$ cheesy polyomino for every $m \in \mathbb{N}$, and the polyomino on the right is not a level $m$ cheesy polyomino for any $m \in \mathbb{N}$.

\begin{figure}
\begin{center}
\includegraphics[width=60mm]{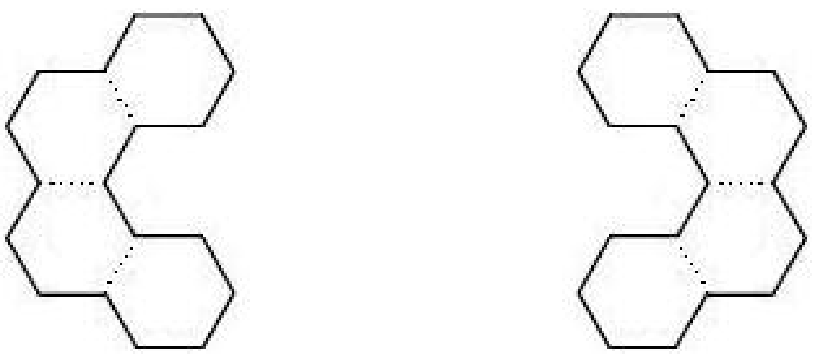}
\caption{The 43rd and 44th four-celled polyominoes.}
\end{center}
\end{figure}

Next: How do the growth constants behave when level tends to infinity? Our database is too small for making precise estimates. Anyway, we know that the growth constant of column-convex polyominoes is $3.863$, while the growth constants of level one, level two and level three cheesy polyominoes are $4.115$, $4.232$ and $4.289$, respectively. Computing the first differences, we get the numbers $4.115-3.863=0.252$, $4.232-4.115=0.117$, and $4.289-4.232=0.057$. Now, the sequence $0.252,\ 0.117,\ 0.057$ is a little similar to a geometric sequence with common ratio $\frac{1}{2}$. Hence, the limit value of the growth constants of cheesy polyominoes might be about $4.232 + 2 \cdot 0.057 = 4.346$.

\section{Conclusion}

This paper is concerned with polyominoes on the hexagonal lattice. For every $m \in \mathbb{N}$, we have defined a set of polyominoes called level $m$ cheesy polyominoes. A polyomino $P$ is a level $m$ cheesy polyomino if the following holds:

\begin{enumerate}
\item $P$ is a rightward-semi-directed polyomino,
\item every column of $P$ has at most two connected components,
\item if a column of $P$ has two connected components, then the gap between the components consists of at most $m$ cells.
\end{enumerate}

Column-convex polyominoes are a subset of level one cheesy polyominoes and, for every $m \in \mathbb{N}$, level $m$ cheesy polyominoes are a subset of level $m+1$ cheesy polyominoes.

For every $m \in \mathbb{N}$, level $m$ cheesy polyominoes have a rational area generating function. We have computed the area generating functions for levels one, two and three. At those three levels, the number of $n$-celled cheesy polyominoes is asymptotically equal to $0.1441\ldots \times 4.1149\ldots ^n$, to $0.1210\ldots \times 4.2318\ldots ^n$, and to $0.1082\ldots \times 4.2886\ldots ^n$, respectively. For comparison, the number of $n$-celled column-convex polyominoes is asymptotically equal to $0.1884\ldots \times 3.8631\ldots ^n$. The number of $n$-celled multi-directed animals behaves asymptotically as $constant \times 4.5878\ldots ^n \:$ \cite{Rechnitzer}. (At present, multi-directed animals are the largest exactly solved class of polyominoes. However, multi-directed animals are not a superset of level one cheesy polyominoes.) The number of all $n$-celled polyominoes behaves asymptotically as $\frac{0.2734\ldots}{n} \times 5.1831\ldots ^n$ \cite{Voege}.

This work could be generalized in several ways. The requirement that polyominoes are rightward-semi-directed can be relaxed and even removed. However, when the above definition is reduced to requirements no. 2 and 3, the area generating function is not a rational function, but a complicated $q$-series. It is also possible not to require rightward-semi-directedness and, at the same time, allow two-component columns to have gaps of all sizes. The only remaining requirement is then the second one, ``every column of $P$ has at most two connected components". The said requirement by itself defines an unsolvable model, but that model can be made solvable by introducing a new requirement. For example, if we forbid runs of two or more consecutive two-component columns, then the area generating function is again a complicated, but computable, $q$-series.

As far as I can see, if requirement no. 2 is replaced by ``every column of $P$ has at most three connected components", the resulting model is still solvable, but if requirement no. 2 is replaced by ``every column of $P$ has at most four connected components", the resulting model is unsolvable.

I think that, already at level one, cheesy polyominoes cannot be enumerated by perimeter. Namely, the perimeter generating function has zero radius of convergence, as can be proved by adapting an argument given by Tony Guttmann in Section 9 of \cite{two}.

\end{document}